\documentclass{amsart}
\usepackage{latexsym}
\usepackage{amsmath, amssymb}
\usepackage[all]{xy}
\newtheorem{dfs}{Definition}[section]

\newtheorem{thms}[dfs]{Theorem}

\newtheorem{blank}[dfs]{}

\newcommand{\dr}{\mathrm{dr}\,}

\begin{document}

\title[Regularity properties for amenable C$^*$-algebras]{Regularity properties 
in the classification program for separable amenable C$^*$-algebras}
\author{George A. Elliott}
\address{Department of Mathematics, University of Toronto, 
Toronto, Ontario, Canada, M5S 2E4}
\email{elliott@math.toronto.edu}
\author{Andrew S. Toms}
\address{Department of Mathematics and Statistics, York University,
4700 Keele St., Toronto, Ontario, Canada, M3J 1P3}
\email{atoms@mathstat.yorku.ca}

\date{\today}

\thanks{This work was partially supported by NSERC}

\begin{abstract}
We report on recent progress in the program to classify separable amenable C$^*$-algebras.
Our emphasis is on the newly apparent role of regularity properties such as finite
decomposition rank, strict comparison of positive elements, and $\mathcal{Z}$-stability, 
and on the importance of the Cuntz semigroup.
We include a brief history of the program's successes since 1989, a more detailed look at the
Villadsen-type algebras which have so dramatically changed the landscape, and a
collection of announcements on the structure and properties of the Cuntz semigroup.

\end{abstract}

\maketitle

\section{Introduction}

Rings of bounded operators on Hilbert space were first studied by
Murray and von Neumann in the 1930s.  These rings, later called von Neumann algebras, 
came to be viewed as a subcategory of a more general category, namely, C$^{*}$-algebras.
(The C$^*$-algebra of compact operators appeared for perhaps the first
time when von Neumann proved the uniqueness of the canonical commutation relations.)
A C$^{*}$-algebra is a Banach algebra $A$ with involution $x \mapsto
x^{*}$ satisfying the C$^{*}$-algebra identity:
\[
||x x^{*}|| = ||x||^{2}, \ \forall x \in A.
\]
Every C$^*$-algebra is isometrically $*$-isomorphic to a norm-closed
sub-$*$-algebra of the $*$-algebra of bounded linear operators on some Hilbert space, 
and so may still be viewed as a ring of operators on a Hilbert space.  

In 1990, the first named author initiated a program to classify
amenable norm-separable C$^{*}$-algebras via $\mathrm{K}$-theoretic
invariants.  The graded and (pre-)ordered group $\mathrm{K}_{0} \oplus
\mathrm{K}_{1}$ was suggested as a first approximation to the correct
invariant, as it had already proved to be complete for both
approximately finite-dimensional (AF) algebras and approximately
circle (A$\mathbb{T}$) algebras of real rank zero (\cite{El1}, \cite{El3}).  It was quickly
realised, however, that more sensitive invariants would be required
if the algebras considered were not sufficiently rich in
projections.  The program was refined, and became concentrated on
proving that Banach algebra $\mathrm{K}$-theory and positive traces
formed a complete invariant for simple separable amenable
C$^{*}$-algebras.  Formulated as such, it enjoyed tremendous 
success throughout the 1990s and early 2000s.

Recent examples based on the pioneering work of Villadsen 
have shown that the classification
program must be further revised.  Two things are now apparent:
the presence of a dichotomy among separable amenable C$^{*}$-algebras
dividing those algebras which are classifiable via $\mathrm{K}$-theory 
and traces from those which will require finer invariants; and the 
possibility---the reality, in some cases---that this dichotomy 
is characterised by one of three potentially equivalent regularity properties for
amenable C$^{*}$-algebras.  (Happily, the vast majority of
our stock-in-trade simple separable amenable C$^{*}$-algebras have
one or more of these properties, including, for instance, those
arising from directed graphs or minimal C$^{*}$-dynamical systems.)

Our plan in this article is to give a brief account of the activity 
in the classification program over the past decade, with 
particular emphasis on the now apparent role of regularity properties.
After reviewing the successes of the program so far, 
we will cover the work of Villadsen on rapid dimension growth 
AH algebras, the examples of R\o rdam and the second named author
which have necessitated the present re-evaluation of the classification 
program, and some recent and sweeping classification results of Winter  
obtained in the presence of the aforementioned regularity
properties.  We will also discuss the possible consequences to the
classification program of including the Cuntz semigroup as part of
the invariant (as a refinement of the $\mathrm{K}_0$ and tracial
invariants).

\section{Preliminaries}

Throughout the sequel $\mathcal{K}$ will denote the C$^*$-algebra of compact operators
on a separable infinite-dimensional Hilbert space $\mathcal{H}$.   
For a C$^*$-algebra $A$, we let $\mathrm{M}_n(A)$ denote the algebra of $n \times n$
matrices with entries from $A$.  The cone of positive elements of $A$ will be denoted by $A_+$.

\subsection{The Elliott invariant and the original conjecture}

The Elliott invariant of a C$^*$-algebra $A$ is the 4-tuple
\begin{equation}\label{ell}
\mathrm{Ell}(A) := \left( (\mathrm{K}_0A,\mathrm{K}_0A^+,\Sigma_A),\mathrm{K}_1A,\mathrm{T}^+A,\rho_A \right),
\end{equation}
where the $\mathrm{K}$-groups are the topological ones, $\mathrm{K}_0A^+$ is the image
of the Murray-von Neumann semigroup $\mathrm{V}(A)$ under the Grothendieck map, $\Sigma_A$ is
the subset of $\mathrm{K}_0A$ corresponding to projections in $A$, $\mathrm{T}^+A$ is
the space of positive tracial linear functionals on $A$, and $\rho_A$ is the natural pairing
of $\mathrm{T}^+A$ and $\mathrm{K}_0A$ given by evaluating a trace at a $\mathrm{K}_0$-class.
The reader is referred to R\o rdam's monograph \cite{R1} for a detailed treatment of this invariant.
In the case of a unital C$^*$-algebra the invariant becomes
\[
\left( (\mathrm{K}_0A,\mathrm{K}_0A^+,[1_A]),\mathrm{K}_1A,\mathrm{T}A,\rho_A \right),
\]
where $[1_A]$ is the $\mathrm{K}_0$-class of the unit, and $\mathrm{T}A$ is the (compact convex)
space of tracial states.  We will concentrate on unital C$^*$-algebras in the sequel
in order to limit technicalities.

The original statement of the classification conjecture for simple unital separable amenable
C$^*$-algebras read as follows:

\begin{blank}
Let $A$ and $B$ be simple unital separable amenable C$^*$-algebras, and suppose that there
exists an isomorphism
\[
\phi:\mathrm{Ell}(A) \to \mathrm{Ell}(B).
\]
It follows that there is a $*$-isomorphsim $\Phi:A \to B$ which induces $\phi$.
\end{blank}

\noindent
It will be convenient to have an abbreviation for the statement above.  Let us call it (EC).

\subsection{Amenability}
We will take the following deep theorem, which combines results of Choi and Effros (\cite{ce}), Connes (\cite{co}), 
Haagerup (\cite{ha}), and Kirchberg (\cite{ki2}), to be our definition of amenability.

\begin{thms}\label{amenabledef}
A C$^*$-algebra $A$ is amenable if and only if it has the following property:  for
each finite subset $\mathcal{G}$ of $A$ and $\epsilon>0$ there are a finite-dimensional
C$^*$-algebra $F$ and completely positive contractions $\phi,\psi$ such that the diagram
\[
\xymatrix{
{A}\ar[rr]^{\mathbf{id}_A}\ar[dr]_{\phi} && {A} \\
&{F}\ar[ur]_{\psi}&
}
\] 
commutes up to $\epsilon$ on $\mathcal{G}$.
\end{thms}

\noindent
The property characterising amenability in Theorem \ref{amenabledef} is known as the
{\it completely positive approximation property}.

Why do we only consider separable and amenable C$^*$-algebras in the classification program?
It stands to reason that if one has no good classification of the weak closures of the GNS 
representations for a class of C$^*$-algebras, then one can hardly expect to classify the
C$^*$-algebras themselves.  These weak closures have separable predual if the C$^*$-algebra is
separable. Connes and Haagerup gave a classification of injective von Neumann algebras with
separable predual (see \cite{co2} and \cite{ha2}), while Choi and Effros established that a C$^*$-algebra is amenable if and only
if the weak closure in each GNS representation is injective (\cite{ce2}).    
Separability and amenability are thus natural conditions which guarantee the existence
of a good classification theory for the weak closures of all GNS representations of a given C$^*$-algebra.
The assumption of amenability has been shown to be necessary by D\u{a}d\u{a}rlat (\cite{da}).

The reader new to the classification program who desires a 
fuller introduction is referred to R\o rdam's excellent monograph \cite{R1}.

\subsection{The Cuntz semigroup}
One of the three regularity properties alluded to in the introduction is
defined in terms of the {\it Cuntz semigroup}, an analogue for positive elements of the
Murray-von Neumann semigroup $\mathrm{V}(A)$.  It is known that this semigroup will
be a vital part of any complete invariant for separable amenable C$^*$-algebras (\cite{To1}).
Given its importance, we present both its original definition, and a modern version
which makes the connection with classical $\mathrm{K}$-theory more transparent.

\begin{dfs}[Cuntz-R\o rdam---see \cite{Cu1} and \cite{R6}]\label{cu1} {\rm
Let $\mathrm{M}_{\infty}(A)$ denote the algebraic limit of the
direct system $(\mathrm{M}_n(A),\phi_n)$, where $\phi_n:\mathrm{M}_n(A) \to \mathrm{M}_{n+1}(A)$
is given by
\[
a \mapsto \left( \begin{array}{cc} a & 0 \\ 0 & 0 \end{array} \right).
\]
Let $\mathrm{M}_{\infty}(A)_+$ (resp. $\mathrm{M}_n(A)_+$)
denote the positive elements in $\mathrm{M}_{\infty}(A)$ (resp. $\mathrm{M}_n(A)$). 
Given $a,b \in \mathrm{M}_{\infty}(A)_+$, we say that $a$ is {\it Cuntz subequivalent} to
$b$ (written $a \precsim b$) if there is a sequence $(v_n)_{n=1}^{\infty}$ of
elements in some $\mathrm{M}_k(A)$ such that
\[
||v_nbv_n^*-a|| \stackrel{n \to \infty}{\longrightarrow} 0.
\]
We say that $a$ and $b$ are {\it Cuntz equivalent} (written $a \sim b$) if
$a \precsim b$ and $b \precsim a$.  This relation is an equivalence relation,
and we write $\langle a \rangle$ for the equivalence class of $a$.  The set
\[
\mathrm{W}(A) := \mathrm{M}_{\infty}(A)_+/ \sim
\] 
becomes a positively ordered Abelian semigroup when equipped with the operation
\[
\langle a \rangle + \langle b \rangle = \langle a \oplus b \rangle
\]
and the partial order
\[
\langle a \rangle \leq \langle b \rangle \Leftrightarrow a \precsim b.
\]
}
\end{dfs}

Definition \ref{cu1} is slightly unnatural, as it fails to consider positive elements
in $A \otimes \mathcal{K}$.  This defect is the result of mimicking the construction
of the Murray-von Neumann semigroup in letter rather than in spirit.  Each projection
in $A \otimes \mathcal{K}$ is equivalent to a projection in some $\mathrm{M}_n(A)$, 
whence $\mathrm{M}_{\infty}(A)$ is large enough to encompass all possible equivalence
classes of projections.  The same is not true, however, of positive elements and Cuntz 
equivalence.  The definition below amounts essentially to replacing $\mathrm{M}_{\infty}(A)$ with 
$A \otimes \mathcal{K}$ in the definition above (this is a theorem), and
also gives a new and very useful characterisation of 
Cuntz subequivalence.  We refer the reader to \cite{La} and \cite{MT} for background 
material on Hilbert C$^{*}$-modules.

Consider $A$ as a (right) Hilbert C$^{*}$-module over itself, and let 
$H_{A}$ denote the countably infinite direct sum of copies of this module.  
There is a one-to-one correspondence between closed countably
generated submodules of $H_{A}$ and hereditary subalgebras of $A
\otimes \mathcal{K}$:  the hereditary subalgebra $B$ corresponds to
the closure of the span of $B H_{A}$.  Since $A$ is separable, $B$ is 
singly hereditarily generated, and it is fairly routine to prove that
any two generators are Cuntz equivalent in the sense of Definition \ref{cu1}.  Thus, passing from
positive elements to Cuntz equivalence classes factors through the
passage from positive elements to the hereditary subalgebras they
generate.  

Let $X$ and $Y$ be closed countably generated submodules of $H_{A}$.
Recall that the compact operators on $H_{A}$ form a C$^{*}$-algebra 
isomorphic to $A \otimes \mathcal{K}$.  Let us say that $X$ is compactly
contained in $Y$ if $X$ is contained in $Y$ and there is 
compact self-adjoint endomorphism of $Y$ which fixes $X$ pointwise.  
Such an endomorphism extends naturally to 
a compact self-adjoint endomorphism of $H_{A}$, and so may be viewed as a self-adjoint
element of $A \otimes \mathcal{K}$.  We write $X \precsim Y$ if each 
closed countably generated compactly contained submodule of $X$ is
isomorphic to such a submodule of $Y$.  

\begin{thms}[Coward-Elliott-Ivanescu, \cite{cei}]\label{cei}
The relation $\precsim$ on Hilbert C$^{*}$ modules defined above, when
viewed as a relation on positive elements in $\mathrm{M}_{\infty}(A)$,
is precisely the relation $\precsim$ of Definition \ref{cu1}.
\end{thms}

Let $[X]$ denote the Cuntz equivalence class of the module $X$.
One may construct a positive partially ordered Abelian semigroup
$\mathcal{C}\mathit{u}(A)$ by endowing the set of countably generated Hilbert
C$^*$-modules over $A$ with the operation
\[
[X] + [Y] := [X \oplus Y]
\]
and the partial order
\[
[X] \leq [Y] \Leftrightarrow X \precsim Y.
\]
The semigroup $\mathcal{C}\mathit{u}(A)$ coincides with $\mathrm{W}(A)$ whenever
$A$ is stable, i.e., $A \otimes \mathcal{K} \cong A$, and has 
some advantages over $\mathrm{W}(A)$ in general.
First, suprema of increasing sequences always exist in $\mathcal{C}\mathit{u}(A)$.  
This leads to the definition of a category including this structure in
which $\mathcal{C}\mathit{u}(A)$ sits as an object, and as a functor into which it is
continuous with respect to inductive limits.  (Definition \ref{cu1}
casts $\mathrm{W}(A)$ as a functor into just the category of partially ordered Abelian
semigroups with zero.  This functor fails to be continuous with respect
to inductive limits.)  Second, it allows one to prove that if $A$ has 
stable rank one, then Cuntz equivalence of positive elements 
simply amounts to isomorphism of the corresponding Hilbert
C$^{*}$-modules.  This has led, via recent work of Brown, Perera, and 
the second named author, to the complete classification of all
countably generated Hilbert C$^{*}$-modules over $A$ via
$\mathrm{K}_{0}$ and traces, and to the classification
of unitary orbits of positive operators in $A \otimes \mathcal{K}$
through recent work of Ciuperca and the first named author (\cite{BPT}, \cite{bt}, \cite{cie}).

Essentially, $\mathrm{W}(A)$ and $\mathcal{C}\mathit{u}(A)$ contain the same information,
but we have chosen to maintain separate notation both to avoid
confusion and because many results in the literature are stated only
for $\mathrm{W}(A)$.  

Cuntz equivalence is often described roughly as the Murray-von Neumann
equivalence of the support projections of positive elements.  
This heuristic is, modulo accounting for projections, precise 
in C$^*$-algebras for which the Elliott invariant is known to be complete (\cite{PT}).
In the stably finite case, one recovers both $\mathrm{K}_0$, the tracial simplex, and
the pairing $\rho$ (see (\ref{ell})) from the Cuntz semigroup, whence the invariant
\[
(\mathcal{C}\mathit{u}(A),\mathrm{K}_1A)
\]
is finer than $\mathrm{Ell}(A)$ in general.  Remarkably, these two invariants determine each other
in a natural way for the largest class of simple separable amenable C$^*$-algebras 
in which (EC) can be expected to hold (\cite{BPT}, \cite{bt}).    

\section{Three regularity properties}\label{reg}

Let us now describe three agreeable properties which a C*-algebra may 
enjoy.  We will see later how virtually all classification theorems
for separable amenable C$^{*}$-algebras via the Elliott invariant
assume, either explicitly or implicitly, one of these properties.

\subsection{Strict comparison}

Our first regularity property---{\it strict comparison}---is one that guarantees, in simple C$^*$-algebras,
that the heuristic view of Cuntz equivalence described at the end of
Section 2 is in fact accurate for
positive elements which are not Cuntz equivalent to projections (see \cite{PT}).
The property is $\mathrm{K}$-theoretic in character.

Let $A$ be a unital C$^*$-algebra, and denote by $\mathrm{QT}(A)$
the space of normalised 2-quasitraces on $A$ (v.\cite[Definition II.1.1]{BH}).
Let $S(\mathrm{W}(A))$ denote the set of additive and order preserving maps $d$ from $\mathrm{W}(A)$ to $\mathbb{R}^+$
having the property that $d(\langle 1_A \rangle) = 1$.
Such maps are called {\it states}.  Given $\tau \in \mathrm{QT}(A)$, one may 
define a map $d_{\tau}:\mathrm{M}_{\infty}(A)_+ \to \mathbb{R}^+$ by
\begin{equation}\label{ldf}
d_{\tau}(a) = \lim_{n \to \infty} \tau(a^{1/n}).
\end{equation}
This map is lower semicontinuous, and depends only on the Cuntz equivalence class
of $a$.  It moreover has the following properties:
\vspace{2mm}
\begin{enumerate}
\item[(i)] if $a \precsim b$, then $d_{\tau}(a) \leq d_{\tau}(b)$;
\item[(ii)] if $a$ and $b$ are orthogonal, then $d_{\tau}(a+b) = d_{\tau}(a)+d_{\tau}(b)$.
\end{enumerate}
\vspace{2mm}
Thus, $d_{\tau}$ defines a state on $\mathrm{W}(A)$.
Such states are called {\it lower semicontinuous dimension functions}, and the set of them 
is denoted by $\mathrm{LDF}(A)$.  If $A$ has the property that $a \precsim b$ whenever $d(a) < d(b)$
for every $d \in \mathrm{LDF}(A)$, then let us say that $A$ has 
{\it strict comparison of positive elements} or simply {\it strict comparison}.

A theorem of Haagerup asserts that every
element of $\mathrm{QT}(A)$ is in fact a trace if $A$ is exact (\cite{ha3}).  All
amenable C$^{*}$-algebras are exact, so we dispense with the
consideration of quasi-traces in the sequel.

\subsection{Finite decomposition rank}
Our second regularity property, introduced by Kirchberg and Winter, is topological
in flavour.  It is based on a noncommutative version of covering dimension
called decomposition rank.

\begin{dfs}[\cite{KW}, Definitions 2.2 and 3.1]
\label{dr-def}
Let $A$ be a separable $C^*$-algebra. 

(i) A completely positive map $\varphi : \bigoplus_{i=1}^s M_{r_i} \to A$ is $n$-decomposable, 
if there is a decomposition $\{1, \ldots, s\} = \coprod_{j=0}^n I_j$ such that the restriction 
of $\varphi$ to $\bigoplus_{i \in I_j} M_{r_i}$ preserves orthogonality for each $j \in \{0, \ldots, n\}$.

(ii) $A$ has decomposition rank $n$, $\dr A = n$, if $n$ is the least integer such that the following 
holds: Given $\{b_1, \ldots, b_m\} \subset A$ and $\epsilon > 0$, there is a completely positive 
approximation $(F, \psi, \varphi)$ for $b_1, \ldots, b_m$ within $\epsilon$ (i.e., $\psi:A \to F$ 
and $\varphi:F \to A$ are completely positive contractions and $\|\varphi \psi (b_i) - b_i\| 
< \epsilon$) such that $\varphi$ is $n$-decomposable. If no such $n$ exists, we write $\dr A = \infty$.  
\end{dfs}

Decomposition rank has good permanence properties.  It behaves well with respect to quotients, 
inductive limits, hereditary subalgebras, unitization and stabilization.
Its topological flavour comes from the fact that it generalises
covering dimension in the commutative case:  if $X$ is a locally compact second countable 
space, then $\dr \mathrm{C}_0(X) = \dim X$.  We refer the reader to \cite{KW} for details.

The regularity property that we are interested in is finite
decomposition rank, expressed by the inequality $\dr < \infty$.
This can only occur in a stably finite C$^{*}$-algebra.  

\subsection{$\mathcal{Z}$-stability}

The Jiang-Su algebra $\mathcal{Z}$ is a simple separable amenable
and infinite-dimensional C$^{*}$-algebra with the same Elliott
invariant as $\mathbb{C}$ (\cite{JS1}).  We say that a second algebra $A$ is
$\mathcal{Z}$-stable if $A \otimes \mathcal{Z} \cong A$.  
$\mathcal{Z}$-stability is our third regularity property.
It is very robust with respect to common constructions (see \cite{TW1}).  

The next theorem shows $\mathcal{Z}$-stability
to be highly relevant to the classification program.  Recall that a 
pre-ordered Abelian group $(G,G^+)$ is said to be {\it weakly unperforated}
if $nx \in G^+ \backslash \{0\}$ implies $x \in G^+$ for any $x \in G$ and
$n \in \mathbb{N}$.
\begin{thms}[Gong-Jiang-Su, \cite{gjs}]\label{gjs1}
Let $A$ be a simple unital C$^*$-algebra with weakly unperforated $\mathrm{K}_0$-group.
Then,
\[
\mathrm{Ell}(A) \cong \mathrm{Ell}(A \otimes \mathcal{Z}).
\]
\end{thms}
\noindent
Thus, modulo a mild restriction on $\mathrm{K}_0$, the completeness of 
$\mathrm{Ell}(\bullet)$ in the simple unital case of the classification
program would imply $\mathcal{Z}$-stability.  Remarkably, there exist 
algebras satisfying the hypotheses of the above theorem which are {\it 
not} $\mathcal{Z}$-stable (\cite{R3}, \cite{To1}, \cite{To2}).  

\subsection{Relationships}

In general, no two of the regularity properties above are equivalent.  The most 
important general result connecting them is the following theorem of M. R\o rdam (\cite{R4}):

\begin{thms}
Let $A$ be a simple, unital, exact, finite, and $\mathcal{Z}$-stable C$^*$-algebra.
Then, $A$ has strict comparison of positive elements.  
\end{thms}

\noindent
We shall see in the sequel that for a substantial class of simple, separable, amenable, and stably
finite C$^*$-algebras, all three of our regularity properties are equivalent.  Moreover, the algebras
in this class which do satisfy these three properties also satisfy (EC).  
There is good reason to believe that the equivalence of these three properties will hold
in much greater generality, at least in the stably finite case;  in the general case, 
strict comparison and $\mathcal{Z}$-stability may well prove to be equivalent characterisations of those
simple, unital, separable, and amenable C$^*$-algebras which satisfy (EC). 

\section{A brief history}

We will now take a short tour of the classification program's
biggest successes, and also the fascinating algebras of Villadsen.  
We have two goals in mind:  to edify the reader unfamiliar 
with the classification program, and to demonstrate that the regularity properties of Section
\ref{reg} pervade the known confirmations of (EC).  
This is a new point of view, for when these results were originally proved, there
was no reason to think that anything more than simplicity,
separability, and amenability would be required to complete the
classification program.

We have divided our review of known classification results into three
broad categories according to the types of algebras covered:  purely
infinite algebras, and two formally different types of stably finite algebras.
It is beyond the scope of this article to provide and exhaustive list of
of known classification results, much less demonstrate their connections to 
our regularity properties.  We will thus choose, from each of the three
categories above, the classification theorem with the broadest scope, and 
indicate how the algebras it covers satisfy at least one of our regularity properties.

\subsection{Purely infinite simple algebras}
We first consider a case where the theory is summarised with one 
beautiful result.  Recall that a simple
separable amenable C$^{*}$-algebra is purely infinite if every
hereditary subalgebra contains an infinite projection (a
projection is infinite if it is equivalent, in the sense of Murray and von
Neumann, to a proper subprojection of itself---otherwise the
projection is finite).

\begin{thms}[Kirchberg-Phillips, 1995, \cite{K} and \cite{P1}]\label{pi}
Let $A$ and $B$ be separable amenable purely infinite
simple C$^{*}$-algebras which satisfy the Universal Coefficient Theorem.
If there is an isomorphism
\[
\phi:\mathrm{Ell}(A) \to \mathrm{Ell}(B),
\]
then there is a $*$-isomorphism $\Phi:A \to B$.
\end{thms}

In the theorem above, the Elliott invariant is somewhat simplified.
The hypotheses on $A$ and $B$ guarantee that they are traceless, and
that the order structure on $\mathrm{K}_{0}$ is irrelevant.  Thus, the
invariant is simply  the graded group $\mathrm{K}_{0} \oplus
\mathrm{K}_{1}$, along with the $\mathrm{K}_{0}$-class of the unit if
it exists.  The assumption of the Universal Coefficient Theorem (UCT) 
is required in order to deduce the theorem from a result which is
formally more general:  $A$ and $B$ as in the theorem are
$*$-isomorphic if and only if they are $\mathrm{KK}$-equivalent.  The 
question of whether every amenable C$^{*}$-algebra satisfies the UCT
is open.

Which of our three regularity properties are present here?  As noted
earlier, finite decomposition rank is out of the question.  The
algebras we are considering are traceless, and so the definition of strict
comparison reduces to the following statement:  for any two non-zero
positive elements $a,b \in A$, we have $a \precsim b$.  This, in turn,
is often taken as the very definition of pure infiniteness, and can be
shown to be equivalent to the definition preceding Theorem \ref{pi}
without much difficulty.  Strict comparison is thus satisfied in a
slightly vacuous way.  As it turns out, $A$ and $B$ are also
$\mathcal{Z}$-stable, although this is less obvious.  One first proves
that $A$ and $B$ are approximately divisible (again, this does not
require Theorem \ref{pi}), and then uses the fact, due to Winter and
the second named author, that any separable and approximately
divisible C$^{*}$-algebra is $\mathcal{Z}$-stable (\cite{TW2}).

\subsection{The stably finite case, I:  inductive limits}
We now move on to the case of stably finite C$^{*}$-algebras, i.e.,
those algebras $A$ such that that every projection in the (unitization
of) each matrix algebra $\mathrm{M}_{n}(A)$ is finite.  (The question 
of whether a simple amenable C$^{*}$-algebra must always be purely
infinite or stably finite was recently settled negatively by R\o rdam.
We will address his example again later.)  Many of the
classification results in this setting apply to classes of
C$^{*}$-algebras which can be realised as inductive limits of certain
building block algebras.  The original classification result for
stably finite algebras is due to Glimm.  Recall that a C$^{*}$-algebra
$A$ is uniformly hyperfinite (UHF) if it is the limit of an inductive 
sequence
\[
\mathrm{M}_{n_{1}} \stackrel{\phi_{1}}{\longrightarrow}
\mathrm{M}_{n_{2}} \stackrel{\phi_{2}}{\longrightarrow}
\mathrm{M}_{n_{3}} \stackrel{\phi_{3}}{\longrightarrow} \cdots,
\]
where each $\phi_{i}$ is a unital $*$-homomorphism.  We will state his
result here as a confirmation of the Elliott conjecture, but note that
it predates both the classification program and the realisation that
$\mathrm{K}$-theory is the essential invariant.

\begin{thms}[Glimm, 1960, \cite{Gl}]\label{uhf}
Let $A$ and $B$ be UHF algebras, and suppose that there is an
isomorphsim
\[
\phi:\mathrm{Ell}(A) \to \mathrm{Ell}(B).
\]
It follows that there is a $*$-isomorphism $\Phi: A \to B$ which induces $\phi$.
\end{thms}

Again, the invariant is dramatically simplified here.  Only the
ordered $\mathrm{K}_{0}$-group is non-trivial.  The strategy of
Glimm's proof (which did not use $\mathrm{K}$-theory explicitly)
was to ``intertwine'' two inductive sequences
$(\mathrm{M}_{n_{i}}, \phi_{i})$ and $(\mathrm{M}_{m_{i}}, \psi_{i})$,
i.e., to find sequences of $*$-homomorphisms $\eta_{i}$ and
$\gamma_{i}$ making the diagram
\[
\xymatrix{
{\mathrm{M}_{n_{1}}}\ar[r]^{\phi_{1}}\ar[d]^{\gamma_{1}} & 
{\mathrm{M}_{n_{2}}}\ar[r]^{\phi_{2}}\ar[d]^{\gamma_{2}} &
{\mathrm{M}_{n_{3}}}\ar[r]^{\phi_{3}}\ar[d]^{\gamma_{3}} & \cdots \\
{\mathrm{M}_{m_{1}}}\ar[r]^{\psi_{1}}\ar[ur]^{\eta_{1}} &
{\mathrm{M}_{m_{2}}}\ar[r]^{\psi_{2}}\ar[ur]^{\eta_{2}} &
{\mathrm{M}_{m_{3}}}\ar[r]^{\psi_{3}}\ar[ur]^{\eta_{3}} & \cdots
}
\]
commute.  One then gets an isomorphism between the limit algebras by
extending the obvious morphism between the inductive sequences by
continuity.  

The intertwining argument above can be pushed surprisingly far.  
One replaces the inductive sequences above with more general inductive
sequences $(A_{i}, \phi_{i})$ and $(B_{i}, \psi_{i})$, where the
$A_{i}$ and $B_{i}$ are drawn from a specified class (matrix algebras 
over circles, for instance), and seeks maps $\eta_{i}$ and
$\gamma_{i}$ as before.  Usually, it is not possible to find
$\eta_{i}$ and $\gamma_{i}$ making the diagram commute, but
approximate commutativity on ever larger finite sets can be arranged
for, and this suffices for the existence of an isomorphism between the
limit algebras.  This generalised intertwining is known as the Elliott
Intertwining Argument.  

The most important classification theorem for inductive limits covers 
the so-called approximately homogeneous (AH) algebras.  An AH algebra 
$A$
is the limit of an inductive sequence $(A_{i},\phi_{i})$, where each
$A_{i}$ is {\it semi-homogeneous}:  
\[
A_i = \bigoplus_{j=1}^{n_{i}} p_{i,j}(\mathrm{C}(X_{i,j}) \otimes \mathcal{K})p_{i,j}
\]
for some natural number $n_{i}$, compact metric spaces $X_{i,j}$,
and projections $p_{i,j} \in \mathrm{C}(X_{i,j}) \otimes \mathcal{K}$.
We refer to the sequence $(A_{i},\phi_{i})$ as a {\it decomposition}
for $A$;  such decompositions are not unique.  All AH algebras are
separable and amenable.  

Let $A$ be a simple unital AH algebra.  Let us say that $A$ has slow
dimension growth if it has a decomposition $(A_{i},\phi_{i})$
satisfying
\[
\limsup_{i \to \infty} \sup \left\{
\frac{\mathrm{dim}(X_{i,1})}{\mathrm{rank}(p_{i,1})},\ldots,
\frac{\mathrm{dim}(X_{i,n_{i}})}{\mathrm{rank}(p_{i,n_{i}})}
\right\} = 0.
\]
Let us say that $A$ has very slow dimension growth if it has a
decomposition satisfying the (formally) stronger condition that
\[
\limsup_{i \to \infty} \sup \left\{
\frac{\mathrm{dim}(X_{i,1})^{3}}{\mathrm{rank}(p_{i,1})},\ldots,
\frac{\mathrm{dim}(X_{i,n_{i}})^{3}}{\mathrm{rank}(p_{i,n_{i}})}
\right\} = 0.
\]
Finally, let us say that $A$ has bounded dimension if there is
a constant $M>0$ and a decomposition of $A$ satisfying
\[
\sup_{{i,l}} \{ \mathrm{dim}(X_{i,l}) \} \leq M.
\]

\begin{thms}[Elliott-Gong and D{\u a}d{\u a}rlat, \cite{EG} and \cite{D}]\label{sdgrr0}
(EC) holds among simple unital AH algebras with slow dimension growth and real rank
zero. 
\end{thms}

\begin{thms}[Elliott-Gong-Li and Gong, \cite{EGL} and \cite{G}]\label{vsdg}
(EC) holds among simple unital AH algebras with very slow dimension growth.
\end{thms}

All three of our regularity properties hold for the algebras of
Theorems \ref{sdgrr0} and \ref{vsdg}, but some are easier to establish
than others.  Let us first point out that an algebra from either class
has stable rank one and weakly unperforated $\mathrm{K}_{0}$-group (cf.\cite{BDR}), 
and that these facts predate Theorems
\ref{sdgrr0} and \ref{vsdg}.  A simple unital C$^{*}$-algebra of real 
rank zero and stable rank one has strict comparison if and only if its
$\mathrm{K}_{0}$-group is weakly unperforated (cf.\cite{pe2}),
whence strict comparison holds for the algebras covered by Theorem 
\ref{sdgrr0}.  A recent result of the second named author shows that
strict comparison holds for any simple unital AH algebra with slow
dimension growth (\cite{To6}), and this result is independent of the
classification theorems above.  Thus, strict comparison holds for the
algebras of Theorems \ref{sdgrr0} and \ref{vsdg}, and the proof of 
this fact, while not easy, is at least much less complicated than the
proofs of the classification theorems themselves.  Establishing finite
decomposition rank requires the full force of the classification
theorems:  a consequence of both theorems is that the algebras they 
cover are all in fact simple unital AH algebras of bounded dimension, 
and such algebras have finite decomposition rank by 
\cite[Corollary 3.12 and 3.3 (ii)]{KW}.  Proving $\mathcal{Z}$-stability is also
an application of Theorems \ref{sdgrr0} and \ref{vsdg}:  one may use the said
theorems to prove that the algebras in question are approximately
divisible (\cite{EGL2}), and this entails $\mathcal{Z}$-stability for
separable C$^{*}$-algebras (\cite{TW2}).

Why all the interest in inductive limits?  Initially at least, 
it was surprising to find that any classification of C$^{*}$-algebras by
$\mathrm{K}$-theory was possible, and the earliest theorems to this
effect covered inductive limits (see, for instance, the first named author's
classification of AF algebras and A$\mathbb{T}$-algebras of real rank 
zero \cite{}).  But it was the realisation by Evans and the first named author 
that a very natural class of C$^{*}$-algebras arising from dynamical
systems---the irrational rotation algebras---
were in fact inductive limits of elementary building
blocks that began the drive to classify inductive limits of all
stripes (\cite{EE}).  This theorem of Elliott and Evans has recently been
generalised in sweeping fashion by Lin and Phillips, who prove that
virtually every C$^{*}$-dynamical system giving rise to a simple
algebra is an inductive limit of fairly tractable building blocks.
This result continues to provide strong motivation for the study of 
inductive limit algebras.

\subsection{The stably finite case, II:  tracial approximation}

Natural examples of separable amenable C$^{*}$-algebras are rarely
equipped with obvious and useful inductive limit decompositions.  
Even the aforementioned theorem of Lin and Phillips, which gives an
inductive limit decomposition for each minimal C$^{*}$-dynamical
system, does not produce inductive sequences covered by 
existing classification theorems.  It is thus desirable to have
theorems confirming the Elliott conjecture under hypotheses that are
(reasonably) straightforward to verify for algebras not given as inductive limits.

Lin in \cite{Li3} introduced the concept of tracial topological rank
for C$^{*}$-algebras.  His definition, in spirit if not in letter, is
this:  a unital simple tracial C$^{*}$-algebra
$A$ has tracial topological rank at most $n \in \mathbb{N}$ if for any finite
set $\mathcal{F} \subseteq A$, tolerance $\epsilon>0$, and positive element $a \in A$ there exist
unital subalgebras $B$ and $C$ of $A$ such that
\begin{enumerate}
\item[(i)] $\mathbf{1}_A = \mathbf{1}_B \oplus \mathbf{1}_C$,
\item[(ii)] $\mathcal{F}$ is almost (to within $\epsilon$) contained in $B \oplus C$,
\item[(iii)] $C$ is isomorphic to $F \otimes \mathrm{C}(X)$, where $\mathrm{dim}(X) \leq n$ and
$F$ is finite-dimensional, and
\item[(iv)] $\mathbf{1}_B$ is dominated, in the sense of Cuntz subequivalence, by $a$.
\end{enumerate}  
One denotes by $\mathrm{TR}(A)$ the least integer $n$ for which $A$ satisfies the definition
above;  this is the tracial topological rank, or simply the tracial rank, of $A$.

The most important value of the tracial rank is zero.
Lin proved that simple unital separable amenable C$^*$-algebras of tracial rank zero
satisfy the Elliott conjecture, modulo the ever present UCT assumption (\cite{Li5}).  The great
advantage of this result is that its hypotheses can be verified for a wide variety
of C$^*$-dynamical systems and all simple non-commutative tori, without ever having to prove that the latter have tractable
inductive limit decompositions (see \cite{ph}, for instance).  Indeed, the existence of such decompositions is a
consequence of Lin's theorem!  One can also verify the hypotheses of Lin's classification 
theorem for many real rank zero C$^*$-algebras with unique trace (\cite{B}), always with
the assumption, indirectly, of strict comparison. 

Simple unital C$^*$-algebras of tracial rank zero can be 
shown to have stable rank one and weakly
unperforated $\mathrm{K}_0$-groups, whence they have strict comparison of positive elements
by a theorem of Perera (\cite{pe2}).  (There is a classification theorem for algebras of 
tracial rank one (\cite{Li2}), but this has been somewhat less useful---it is difficult to 
verify tracial rank one in natural examples.  Also, Niu has recently proved a classification
theorem for some C$^*$-algebras which are approximated in trace by certain subalgebras
of $\mathrm{M}_n \otimes \mathrm{C}[0,1]$ (\cite{ni}).)  

And what of our regularity properties? 
Lin proved in \cite{Li3} that every unital simple C$^{*}$-algebra of 
tracial rank zero has stable rank one and weakly unperforated
$\mathrm{K}_{0}$-group.  These facts, by the results reviewed at the end of
the preceding subsection, entail strict comparison, and are not nearly
so difficult to prove as the tracial rank zero classification theorem.  
In a further
analogy with the case of AH algebras, finite decomposition rank and
$\mathcal{Z}$-stability can only be verified by applying Lin's
classification theorem---a consequence of this theorem is that the
algebras it covers are in fact AH algebras of bounded dimension!

\subsection{Villadsen's algebras}

Until the mid 1990s we had no examples of simple separable amenable C$^*$-algebras where one of our
regularity properties failed.  To be fair, two of our regularity properties had not yet even been defined, 
and strict comparison was seen as a technical version of the more attractive Second Fundamental Comparability
Question for projections (this last condition, abbreviated FCQ2, asks for strict comparison for projections only).  
This all changed when Villadsen produced a simple separable amenable and stably finite C$^*$-algebra which did not have FCQ2,
answering a long-standing question of Blackadar (\cite{V1}).  The techniques introduced by Villadsen were subsequently used by 
him and others to answer many open questions in the theory of nuclear C$^*$-algebras including the following:
\begin{enumerate}
\item[(i)] Does there exist a simple separable amenable C$^*$-algebra containing a finite and an infinite projection?
(Solved affirmatively by R\o rdam in \cite{R3}.)
\item[(ii)] Does there exist a simple and stably finite C$^*$-algebra with non-minimal stable rank?
(Solved affirmatively by Villadsen in \cite{V2}.)
\item[(iii)] Is stability a stable property for simple C$^*$-algebras?
(Solved negatively by R\o rdam in \cite{R5}.)
\item[(iv)] Does a simple and stably finite C$^*$-algebra with cancellation of projections necessarily
have stable rank one? (Solved negatively by the second named author in \cite{To5}.)
\end{enumerate}
Of the results above, (i) was (and is) the most significant.  In addition to showing that simple separable amenable C$^*$-algebras do
not have a factor-like type classification, R\o rdam's example demonstrated that the Elliott invariant as it
stood could not be complete in the simple case.  This and other examples due to the second named author have necessitated
a revision of the classification program.  It is to the nature of this revision 
that we now turn.

\section{The way(s) forward}


\subsection{New assumptions}\label{assump}

(EC) does not hold in general, and this justifies new assumptions
in efforts to confirm it.  In particular, one may assume any combination
of our three regularity properties.  We will comment on the aptness of
these new assumptions in the next subsection.  For now we observe that, 
from a certain point of view, we have been making these assumptions all along.  Existing
classification theorems for C$^{*}$-algebras of real rank zero are
accompanied by the crucial assumptions of stable rank one and weakly
unperforated $\mathrm{K}$-theory;  as has already been pointed out,
unperforated $\mathrm{K}$-theory can be replaced with strict comparison 
in this setting.  

How much further
can one get by assuming the (formally) stronger condition of $\mathcal{Z}$-stability?
What role does finite decomposition rank play?  As it turns out, these 
two properties both alone and together produce interesting results.
Let $\mathcal{RR}0$ denote the class of simple unital separable
amenable C$^{*}$-algebras of real rank zero.  The
following subclasses of $\mathcal{RR}0$ satisfy (EC):
\begin{enumerate}
\item[(i)] algebras which satisfy the UCT, have finite decomposition
rank, and have tracial simplex with compact and zero-dimensional extreme 
boundary;
\item[(ii)] $\mathcal{Z}$-stable algebras which satisfy the UCT and 
are approximated locally by subalgebras of finite decomposition rank.
\end{enumerate}
\noindent
These results, due to Winter (\cite{Wi1}, \cite{Wi2}), showcase 
the power of our regularity properties:  included in the algebras
covered by (ii) are all simple separable unital $\mathcal{Z}$-stable
ASH (approximately subhomogeneous) algebras of real rank zero.  

Another advantage to the assumptions of $\mathcal{Z}$-stability and
strict comparison is that they allow one to recover extremely fine
isomorphism invariants for C$^*$-algebras from the Elliott invariant
alone.  (This recovery is not possible in general.)  We will be able to give 
precise meaning to this comment below, but first require a further
dicussion of Cuntz semigroup.   

\subsection{New invariants}

A natural reaction to an incomplete invariant is to enlarge it:  
include whatever information was used to prove incompleteness.  
This is not always a good idea.  It is possible that one's distinguishing information is 
{\it ad hoc}, and unlikely to yield a complete invariant.  Worse, one may throw so
much new information into the invariant that the impact of its potential completeness is
severely diminished.  The revision of an invariant is a delicate business.  In this light,
not all counterexamples are equal.  

R\o rdam's finite-and-infinite-projection example is distinguished from a simple and purely
infinite algebra with the same $\mathrm{K}$-theory by the obvious fact that the latter
contains no finite projections.  The natural invariant which captures this difference
is the semigroup of Murray-von Neumann equivalence classes of projections in matrices 
over an algebra $A$, denoted by $\mathrm{V}(A)$.  After the appearance of R\o rdam's example, 
the second named author produced a pair of simple, separable, amenable, and stably finite
C$^*$-algebras which agreed on the Elliott invariant, but were not isomorphic.  In this 
case the distinguishing invariant was Rieffel's stable rank.  It was later discovered that
these algebras could not be distinguished by their Murray-von Neumann semigroups, but it was
not yet clear which data were missing from the Elliott invariant.  More dramatic
examples were needed, ones which agreed on most candidates for enlarging the invariant, and
pointed the way to the ``missing information''.

In \cite{To2}, the second named author constructed a pair of simple unital AH algebras which, while
non-isomorphic, agreed on a wide swath of invariants including the Elliott invariant, all 
continuous (with respect to inductive sequences) and homotopy invariant functors from the category of
C$^*$-algebras (a class which includes the Murray-von Neumann semigroup), 
the real and stable ranks, and, as was shown later in \cite{}, stable isomorphism invariants 
(those invariants which are insensitive to tensoring with a matrix algebra or passing to 
a hereditary subalgebra).  It was thus reasonable
to expect that the distinguishing invariant in this example---the Cuntz semigroup---might be a good candidate
for enlarging the invariant.  At least, it was an object which after years of being used sparingly
as a means to other ends, merited study for its own sake.  

Let us collect some evidence supporting the addition of the Cuntz semigroup to the usual
Elliott invarariant.  First, in the biggest class of algebras where (EC) can
be expected to hold---$\mathcal{Z}$-stable algebras, as shown by Theorem \ref{gjs1}---it is
not an addition at all!  Recent work of Brown, Perera, and the second named author shows that for
a simple unital separable amenable C$^*$-algebra which absorbs $\mathcal{Z}$ tensorially, there is
a functor which recovers the Cuntz semigroup from the Elliott invariant (\cite{BPT}, \cite{PT}).  This functorial recovery
also holds for simple unital AH algebras of slow dimension growth, a class for which $\mathcal{Z}$-stability
is not known and yet confirmation of (EC) is expected.
(It should be noted that the computation of the Cuntz semigroup for a simple approximately interval (AI)
algebra was essentially carried out by Ivanescu and the first named author in \cite{ei},
although one does require \cite[Corollary 4]{cei} to see that the computation is complete.)
Second, the Cuntz semigroup unifies the counterexamples of R\o rdam and
the second named author.  One can show that the examples of \cite{R1}, \cite{To1}, and \cite{To2} all
consist of pairs of algebras with different Cuntz semigroups;  there are no counterexamples to
the conjecture that simple separable amenable C$^*$-algebras will be classified up to $*$-isomorphism
by the Elliott invariant and the Cuntz semigroup.  Third, the Cuntz semigroup provides a 
bridge to the classification of non-simple algebras.  Ciuperca and the
first named author have recently proved that AI algebras---limits of 
inductive sequences of algebras of the form
\[
\bigoplus_{i=1}^n \mathrm{M}_{m_i}(\mathrm{C}[0,1])
\]
--- are classified up to isomorphism by their Cuntz semigroups.  This is
accomplished by proving that the approximate unitary equivalence
classes of positive operators in the unitization of a stable
C$^{*}$-algebra of stable rank one are determined by the Cuntz
semigroup of the algebra, and then appealing to a theorem of Thomsen (\cite{Th}).
(These approximate unitary equivalence classes of positive operators
can be endowed with the structure of a topological partially ordered
semigroup with functional calculus.  This invariant, known as Thomsen's semigroup, is
recovered functorially from the Cuntz semigroup for separable algebras
of stable rank one, and so from the Elliott invariant in algebras
which are moreover simple,
unital, exact, finite, and $\mathcal{Z}$-stable by the
results of \cite{BPT}.  This new semigroup is the fine invariant
alluded to at the end of subsection \ref{assump}.)

There is one last reason to suspect a deep connection between the classification program and the Cuntz semigroup.
Let us first recall a theorem of Kirchberg, which is germane to the classification of purely infinite
C$^*$-algebras (cf. Theorem \ref{pi}).
\begin{thms}[Kirchberg, c.~1994; see \cite{kp}]
Let $A$ be a separable amenable C$^*$-algebra.  The following two properties are equivalent:
\begin{enumerate}
\item[(i)] $A$ is purely infinite;
\item[(ii)] $A \otimes \mathcal{O}_{\infty} \cong A$.
\end{enumerate} 
\end{thms}
\noindent
A consequence of Kirchberg's theorem is that among simple separable amenable C$^*$-algebras
which merely contain an infinite projection, there is a two-fold characterisation of the
(proper) subclass which satisfies the original form of the Elliott conjecture (modulo UCT).  If one assumes {\it a priori} that $A$
is simple and unital with no tracial state, then a theorem of R\o rdam (see \cite{R4}) shows that property (ii) above ---
known as $\mathcal{O}_{\infty}$-stability---is equivalent to $\mathcal{Z}$-stability.
Under these same hypotheses, property (i) is equivalent to the statement that $A$ has strict
comparison.  Kirchberg's theorem can thus be rephrased as follows in the simple unital case:
\begin{thms}
Let $A$ be a simple separable unital amenable C$^*$-algebra without a tracial
state.  The following two properties are equivalent:
\begin{enumerate}
\item[(i)] $A$ has strict comparison;
\item[(ii)] $A \otimes \mathcal{Z} \cong A$.
\end{enumerate} 
\end{thms}
\noindent
The properties (i) and (ii) in the theorem above make perfect sense in the presence of a trace.
We moreover have that (ii) implies (i) even in the presence of traces (this is due to R{\o}rdam---see \cite{R4}).
It therefore makes sense to ask whether the theorem might be true without the tracelessness hypothesis.
Remarkably, this appears to be the case.  Winter and the second named author have proved
that for a substantial class of stably finite C$^*$-algebras, strict comparison and $\mathcal{Z}$-stability
are equivalent, and that these properties moreover characterise the (proper) subclass which
satisfies (EC) (\cite{TW3}).  In other words, Kirchberg's theorem
is quite possibly a special case of a more general result, one which will give a unified two-fold 
characterisation of those simple separable amenable C$^*$-algebras which satisfy the original form
of the Elliott conjecture.

It is too soon to know whether the Cuntz semigroup together with Elliott invariant will
suffice for the classification of simple separable amenable C$^*$-algebras, or indeed,
whether such a broad classification can be hoped for at all.  But there is already cause for optimism.
Zhuang Niu has recently obtained some results on lifting maps at the level of the Cuntz semigroup
to $*$-homomorphisms.  This type of lifting result is a key ingredient in proving
classification theorems of all stripes.  His results suggest the algebras
of \cite{To2} as the appropriate starting point for any effort to establish the Cuntz
semigroup as a complete isomorphism invariant, at least in the absence of $\mathrm{K}_1$.    

We close our survey with a few questions for the future, both near and far. 
\begin{enumerate}
\item[(i)] When do natural examples of simple separable amenable C$^*$-algebras satisfy
one or more of the regularity properties of Section 3?  In particular, do simple unital 
inductive limits of recursive subhomogeneous algebras have strict comparison whenever they
have strict slow dimension growth?
\item[(ii)] Can the classification of positive operators up to approximate unitary
equivalence via the Cuntz semigroup in algebras of stable rank one be extended to
normal elements, provided that one accounts for $\mathrm{K}_1$?
\item[(iii)] Let $A$ be a simple, unital, separable, and amenable C$^*$-algebra with
strict comparison of positive elements.  Is $A$ $\mathcal{Z}$-stable?  Less ambitiously,
does $A$ have stable rank one whenever it is stably finite?
\item[(iv)] Can one use Thomsen's semigroup to prove new classification theorems?  (The attraction here
is that Thomsen's semigroup is already implicit in the Elliott invariant for many classes of C$^*$-algebras.) 
\end{enumerate}

\end{document}